Баяк И.В.

# Об одной математической конструкции многомерного космоса

***Резюме:*** Установлено, что геометрия наблюдаемого евклидова пространства, алгебры Паули и Дирака, а также группы внутренних симметрий стандартной модели и модели так называемого большого объединения калибровочных взаимодействий логически следуют из геометрических построений 8-мерного вакуума и 9-мерного космоса.

Прежде всего, покажем, что классические алгебры Клиффорда являются групповыми алгебрами подгрупп циклических стрелочных подстановок. Пусть $I = \{1,..,n\}$, тогда имеем простое отражение $\mathrm{O}: I \to I : i \to i^\circ$: $\{1,..,i,..,n\} \to \{n_{(1)},...,i^\circ_{(i)},...,1_{(n)}\}$. Пусть также здесь и всюду далее $I = \{1,..,2^m\}$, $J = \{1,..,m\}$, тогда существует такое разбиение семейства $I$ на подсемейства, что $\forall k \in J$ имеем $I^k \equiv \{I^k_j\}_{(2^{k-1})}$, где $\bigcap_{(2^{k-1})} I^k_j = \varnothing$, $\bigcup_{(2^{k-1})} I^k_j = I$, $\forall j \in \{1,..,2^{k-1}\}$ $\left(card\left(I^k_j\right) = 2^{m-k+1}\right)$, и определена индукция, так что каждое $I^k_{j+1}$ заполняется последовательной выборкой из $I$ вслед за $I_j$, причем $I^k_1 = \{1,..,2^{m-k+1}\}$. В свою очередь, для такого $I^k$ существует отображение $k$–отражения, а именно $k \circ : I \to \left\{\mathrm{O}\left(I^k_j\right)\right\}_{(2^{k-1})}$.

Далее пусть всякая подстановка $p$ из группы подстановок $P$ задается множеством биективных переходных пар $(i,j): i \to j$, т.е. $p = \{(i,j) | \forall i \neq i'(j \neq j'), \forall j \neq j'(i \neq i')\}_{i,j \in I}$. Выделим из этой группы множество циклических подстановок $b = \{(i,j)(j,i) | \forall i \neq i'(j \neq j')\}_I$; тем самым $b^2 = e$, где $e$ - тождественная подстановка, а $\{b\}$ - множество циклических подстановок и только. Если же $bp \equiv \{(i,j),(j,i),\ (i^\circ,j^\circ),(j^\circ,i^\circ) | i \neq j, \forall i \neq i'(j \neq j')\}_I$, то $bP = \{bp\}$ - группа циклических подстановок. Вместе с тем, если $bp(k) \equiv \{(i,j) | (i,j) \in bp, j = i^{k\circ}\}_I$, то $\{bp(k)\}_J$ - генераторы группы циклических подстановок, т.е. $bP = \langle \{bp(k)\}_J \rangle$, причем $\forall k_1 \neq k_2$ $(bp(k_1)bp(k_2) = bp(k_2)bp(k_1))$, так что $\forall bp_1 \neq bp_2$ $(bp_1 bp_2 = bp_2 bp_1)$.

Далее пусть всякая подстановка $\bar{p}$ из группы стрелочных подстановок $\bar{P}$, где $P < \bar{P}$, задается множеством знако-переходных пар, т.е. $\bar{p} = \{(i, \pm j) | \forall i \neq i'(j \neq j'), \forall j \neq j'(i \neq i')\}_{i,j \in I}$. Тогда, если $b\bar{p} \equiv \{(i,j),(j,-i),(i^\circ,j^\circ),(j^\circ,-i^\circ) | i \neq j, \forall i \neq i'(j \neq j')\}_I$, то имеем $b\bar{P} = \{b\bar{p}\}$ - группу циклических стрелочных подстановок, где $b\bar{p}^2 = -e$. В свою очередь, если $b\bar{p}(k) \equiv \{(i,\pm j)|$ $(i, \pm j) \in b\overline{p}, j = i^{k\circ}\}_J$, то $\{b\overline{p}(k)\}_J$ - генераторы группы циклических стрелочных подстановок, т.е. $b\overline{P} = \langle \{b\overline{p}(k)\}_J \rangle$, причем $\forall k_1 \neq k_2 \quad (b\overline{p}(k_1)b\overline{p}(k_2) = -b\overline{p}(k_2)b\overline{p}(k_1))$, так что $\forall b\overline{p}_1 \neq b\overline{p}_2 \quad (b\overline{p}_1 b\overline{p}_2 = -b\overline{p}_2 b\overline{p}_1)$.

Таким образом, поскольку свойства генераторов группы циклических стрелочных подстановок соответствуют требованиям, предъявляемым к генераторам классической алгебры Клиффорда, а именно $e_i e_j + e_j e_i = -2\delta_{ij} \cdot 1$, где $i, j \in J$, то групповая алгебра над $R$ с базисом $b\overline{P}$ изоморфна именованной алгебре, т.е. $Rb\overline{P} \approx Cl(2^m)$. Вместе с тем, если для структурных отношений генераторов алгебры Клиффорда потребовать выполнения условия $e_1^2 = 1$, то наряду с классической алгеброй следует рассматривать и алгебру $Rb\hat{P} \approx Cl(2^{m-1}, 2^{m-1})$, где $b\hat{P} = \{b\hat{p}\} = \langle bp(1), \{b\overline{p}(k)\}_{J\setminus 1} \rangle$.

Далее заметим, что $\langle b\overline{p}(k) \rangle \approx P^+(2)$, т.е. циклическая группа, порожденная произвольным генератором циклических стрелочных подстановок, изоморфна подгруппе четных стрелочных подстановок 2-го порядка. В свою очередь, для циклической группы, порожденной генератором $bp(1)$, имеем изоморфизм $\langle bp(1) \rangle \approx P^{\pm}(1,1)$, где $P^{\pm}(1,1)$ - подгруппа четно-составных стрелочных подстановок 2-го порядка. Отсюда немедленно следует, что $b\overline{P} < P^+(2^m)$ а $b\hat{P} < P^{\pm}(2^{m-1}, 2^{m-1})$, и поэтому $Cl(2^m) \subset RP^+(2^m) \approx End\ E(2^m)$ а $Cl(2^{m-1}, 2^{m-1}) \subset RP^{\pm}(2^{m-1}, 2^{m-1}) \approx End\ E(2^{m-1}, 2^{m-1})$, где $E(2^m)$, $E(2^{m-1}, 2^{m-1})$ - соответствующие евклидовы пространства.

Обратимся теперь к алгебрам физического мира, а именно, к алгебрам Паули и Дирака, которые изоморфны алгебрам $Cl(4)$ и $Cl(4,4)$ соответственно, причем первая алгебра Клиффорда действует в 4-мерном евклидовом пространстве, а вторая – в 8-мерном псевдоевклидовом пространстве. Если данному изоморфизму придать определенный физический смысл, тогда необходимо будет коренным образом поменять представление о геометрии физического вакуума. Действительно, обычно геометрию пустого пространства сопоставляют с геометрией наблюдаемого 3-мерного евклидова пространства, но если распространить микроструктуру Вселенной на ее макроструктуру, тогда необходимо дополнительно к этой обычной ассоциации применить в первом приближении геометрию $E(4)$, а во втором приближении – геометрию $E(4,4)$. В свою очередь, геометрия $E(4)$ обусловлена факторизацией $R^4$ в $S^4$, а геометрия $E(4,4)$ - факторизацией $R^8$ в $S^4 \times S^4$, и поэтому постулируем еще одну конструкцию вакуума, а именно, расслоение с базой $E(3)$ и типичным слоем $S^4$ или $S^4 \times S^4$.

Продолжая анализ глобальных физических симметрий, обратимся к группам $SU(3) \times SU(2) \times U(1)$, $SU(5)$, которые описывают идеальные симметрии соответственно стандартной модели и модели так называемого большого объединения калибровочных взаимодействий. При этом в действительности имеет место нарушение идеальных симметрий до группы $SU(3) \times U(1)$. Предположим также, что на самом деле всякая физическая унитарная группа $SU(n)$ должна быть расширена до группы $SO(2n)$, так что геометрическим источником физических симметрий служит орбита действия $SO(2n)$ в $R^{2n}$, т.е. сфера $S^{2n-1}$. В связи с этим нелишне будет напомнить, что овеществление группы $SU(n)$ является подгруппой $SO(2n)$, и что $S^n \approx R^n / ZH_n^+$, причем $End\, S^n \approx End\, E(n) = RP_n^+$ и $Aut\, S^n \approx Aut\, E(n) = SO(n)$, где $ZH_n^+ < Z^n$ а $RP_n^+$ - групповая алгебра над $R$ с базисом $P_n^+$ - подгруппой четных стрелочных подстановок $n$-го порядка. Возвращаясь к геометрическим источникам идеальных физических симметрий, заметим, что $S^5 \oplus S^3 \oplus S^1 < S^9$ (т.к. $R^5 \oplus R^3 \oplus R^1 = R^9$, но $ZH_5^+ \oplus ZH_3^+ \oplus 2Z < ZH_9^+$), и поэтому предметом нашего поиска должно быть такое линейное преобразование факторпространства $S^9$, в результате которого сохраняется симметрия $S^5, S^1$ и нарушается симметрия $S^3$. Однако, поскольку эндоморфизмы евклидовых пространств совпадают с эндоморфизмами сфер (представленных факторпространствами) и $E(9) = E(5) \oplus E(3) \oplus E(1)$, то требуется найти такой эндоморфизм $E(9)$, который сохранит автоморфизмы подпространств $E(5), E(1)$ и нарушит группу автоморфизмов подпространства $E(3)$. Так, например, эндоморфизм подобия с коэффициентами $\alpha, \beta$ и 1, заданный матрицей $diag[\alpha,\alpha,\alpha,\alpha,\alpha,\beta,\beta,\beta,1]$, оставляет инвариантными и $E(5)$ и $E(3)$ и $E(1)$. Если же в качестве линейного преобразования $f$ пространства $E(9)$ принять отображение, заданное матрицей $F$,

$$F \equiv \begin{pmatrix} \alpha & 0 & 0 & 0 & 0 & 0 & 0 & 0 & 0 \\ 0 & \alpha & 0 & 0 & 0 & 0 & 0 & 0 & 0 \\ 0 & 0 & \alpha & 0 & 0 & 0 & 0 & 0 & 0 \\ 0 & 0 & 0 & \alpha & 0 & 0 & 0 & 0 & 0 \\ 0 & 0 & 0 & 0 & \alpha & 0 & 0 & 0 & 0 \\ 0 & 0 & 0 & 0 & 0 & \beta & 0 & 0 & 0 \\ 0 & 0 & 0 & 0 & 0 & 0 & \beta & 0 & 0 \\ 0 & 0 & 0 & 0 & 0 & 0 & 0 & \beta & 0 \\ 0 & 0 & 0 & 0 & 0 & \chi & \chi & \chi & 1 \end{pmatrix},$$

тогда $E(9)$ содержит инвариантные относительно $f$ подпространства $E(5), E(1)$ с ограничением $f(E(5)), f(E(1))$, являющимся эндоморфизмом подобия с коэффициентом 1 и $\alpha$ соответственно, а также потенциально инвариантное (при условии $\chi = 0$) подпространство $E(3)$ с ограничением $f(E(3))$, являющимся эндоморфизмом подобия с коэффициентом $\beta$. Тем самым, настоящий выбор $F$ гарантирует связь стандартной модели теоретической физики с евклидовой геометрией, производной от унитарных симметрий.

Таким образом, математическое соответствие между геометрическими источниками симметрий и унитарными группами коренным образом меняет представление о геометрии физического космоса. Действительно, обычно унитарные симметрии связывают с микрокосмосом, однако если их отнести к космосу в целом, то возникает необходимость рассмотреть конструкцию космоса, основанную на геометрии 9-мерной сферы, а именно, расслоение с базой $E(3)$ и типичным слоем $S^9$.

      В заключение заметим, что наблюдаемое 3-мерное евклидово пространство, которое является базой и вакуумного и космического расслоений, формально может быть факторизовано в 3-сферу и вложено в типичный слой этих расслоений. В свою очередь, вакуумное расслоение вкладывается в космическое, причем $S^4 \times S^4$ глобально минимально в $S^9$, и поэтому можно предположить, что оно эволюционирует к абсолютно минимальному в $S^9$ подпространству $S^8$. В целом, можно также полагать, что вакуумное расслоение - это всего лишь абстрактная формализация 8-мерного слоения с типичным слоем, обертывающим $S^4 \times S^4$, заданного в 9-мерной сфере, подвергнутой вышеуказанному эндоморфизму. Аналогично, можно полагать, что наблюдаемое евклидово пространство – это абстрактная формализация 3-мерного слоения, вложенного в вакуумное слоение и обертывающего 3-сферу в каждом своем слое.